%% file: paper.tex
\documentclass[conference,12pt,onecolumn,draftclsnofoot]{IEEEtran}
\ifCLASSINFOpdf
\else
\fi
\usepackage{amsmath}
\usepackage{amssymb}
\usepackage{amsthm}
\usepackage{bm}
\usepackage[makeroom]{cancel}

\ifCLASSOPTIONcompsoc
  \usepackage[caption=false,font=normalsize,labelfont=sf,textfont=sf]{subfig}
\else
  \usepackage[caption=false,font=footnotesize]{subfig}
\fi

\usepackage[utf8]{inputenc}
\usepackage[T1]{fontenc}

\usepackage{tikz}
\usetikzlibrary{positioning}

\usepackage{xcolor}

\newcommand\numberthis{\addtocounter{equation}{1}\tag{\theequation}}

\newtheorem{definition}{Definition}
\newtheorem{remark}{Remark}

\newcommand{\Real}[0]{\mathbb{R}}

\newcommand{\vetor}[1]{\bm{#1}}
\newcommand{\ones}[0]{\mathbf{1}}

\newcommand{\minimum}[1]{\underline{#1}}
\newcommand{\maximum}[1]{\overline{#1}}

\newcommand{\transpose}[1]{#1^\intercal}
\newcommand{\lagrangian}{\mathcal{L}}

\begin{document}
\IEEEoverridecommandlockouts

\IEEEpubid{978-1-7281-1156-8/19/\$31.00~\copyright{}2019 IEEE \hfill}
%
\title{Financial storage rights for hydroelectricity}

\author{\IEEEauthorblockN{L.~S.~A.~Martins}
\IEEEauthorblockA{IBM Research\\
São Paulo, Brazil 04007--900}
\and
\IEEEauthorblockN{R.~L.~Hochstetler}
\IEEEauthorblockA{Instituto Acende Brasil\\
São Paulo, Brazil 04534--004}}


%


\maketitle

\begin{abstract}
There has recently been growing interest in the development of financial storage rights for energy storage systems as instruments akin to their transmission counterparts as a means to not only distribute congestion rents, but also to mitigate price risks, and even to signal and finance investments in light of increasing penetration of renewable generation and downward pressures on market clearing prices. This work presents a discussion on the applicability of such instruments to hydroelectric power plants as a mechanism to decouple ownership and operation of reservoirs, especially those with small regulating capacity, aiming to improve their valuation and investment risk management in electricity markets. The resulting model takes into consideration reasonable assumptions regarding their nonlinear physics and short-term operation properties, extending the applicability of a strong theoretical framework recently proposed in the literature to a large share of the energy mix.
\end{abstract}


%
\IEEEpeerreviewmaketitle

\input{notation}
\input{introduction}
\input{hydro}
\input{fsr}
\input{discussion}
\input{conclusion}

\bibliographystyle{IEEEtran}
\bibliography{IEEEabrv,references}
%



\end{document}

%% file: notation.tex
\section*{Notation and nomenclature}
Scalars, vectors, and matrices are represented by lowercase, bold lowercase, and uppercase letters, respectively, e.g., $x$, $\vetor{x}$, $X$. The $i$-th component of a vector $\vetor{x}$ is denoted by $x_i$, $X_{i,j}$ denotes the $(i,j)$-th element of a given matrix $X$, and the $j$-th column vector of a matrix $X$ is given by its respective lower case letter as in $\vetor{x}[j]$. All vectors are column vectors, unless noted otherwise, and their transposes are denoted $\transpose{\vetor{x}}$. Vectors of all zeros and ones are denoted $\vetor{0}$ and $\vetor{1}$, respectively.

\begin{description}
	\item[$\lambda$] Locational marginal price.
	\item[$\gamma,\mu,\eta$] Dual variables.
	\item[$\varphi(.)$] Upper reservoir volume-to-height function.
	\item[$\vartheta(.)$] Tailrace discharge-to-height function.
	\item[$\pi$] Hydroelectric productivity.
	\item[$c(.)$] Cost (benefit) function.
	\item[$d$] Line capacity.
	\item[$e$] Energy capacity right.
	\item[$f$] Flow-gate right.
	\item[$G$] Shift factors matrix.
	\item[$h(.)$] Net water head.
	\item[$i,t$] Storage/bus and time indices.
	\item[$p(\cdot)$] Hydroelectric power generation.
	\item[$q$] Water discharge.
	\item[$r$] Financial transmission right.
	\item[$s$] Financial storage right.
	\item[$u$] Power discharge.
	\item[$v$] Water volume stored in upper reservoir.
	\item[$y$] Incremental water inflow to reservoir.
	\item[$z$] Energy storage.
\end{description}

%% file: introduction.tex
\section{Introduction}%
\label{sec:introduction}
\subsection{Motivation and background}%
\label{sec:introduction:motivation}
Most deregulated electric power systems have been restructured around day-ahead and balancing markets. Electricity producers bid price-quantity pairs for each time interval of the following day. If the bid is accepted in the day-ahead auction, the power producer must provide the quantity committed. This arrangement is well-suited for capacity constrained power producers for which it is optimal to produce as long the price is sufficient to cover marginal costs. The same logic, however, does not apply to energy-constrained producers. Although also limited by their installed capacity, the main constraint of a e.g. hydroelectric power plant, is the availability of water stored in its reservoir. Their production optimization is achieved by arbitraging prices across different time periods. Thus, the marginal cost of a hydro plant is given by the inter-temporal opportunity cost of water---its reservoir is depleted when it is dispatched, thus reducing its capacity to meet future demand, possibly implying higher opportunity costs, which are themselves more difficult to assess given the nonlinear productivity of hydroelectric power plants. 

Day-ahead markets are typically not structured in a way that enables energy-constrained producers to adjust supply bids to reflect their current marginal costs. Bids for all time periods of the following day are typically submitted simultaneously with no provision for contingent bids on their reservoir levels. This imposes a daunting task on these producers to determine their optimal bidding strategies~\cite{steeger:2014,aasgard:2018}, since it requires them to base their decisions on projections of market supply and demand dynamics to compute the expected prices for each of the bidding time periods of the next day. This becomes particularly difficult in a market with large share of renewable energy sources (RES), increasing stochastic supply-side variations.

\subsection{Relevant literature}%
\label{sec:introduction:literature}
\IEEEpubidadjcol
Market pricing coupled with financial instruments have been successfully deployed in many power markets as a way to promote optimal operation while providing financial stability to market players. One of the forerunners of this type of arrangement is the sale of financial property rights to scarce transmission capacity~\cite{hogan:1992}. This arrangement was introduced after deregulation and restructuring of electricity markets, most notably in the US, as a hedging instrument for locational marginal price (LMP) risks that arise from network congestion. The buyer of the financial property right receives the network congestion rents, thus providing protection from price differentials that may arise in market transactions across the network, while the sale of these rights provides a steady income stream that helps remunerate investment in transmission~\cite{oren:2013,kristiansen:2006}. These financial property rights have historically taken on two different contractual forms: the point-to-point financial transmission rights (FTR), with which a holder hedges against price differences between two specific nodes; and shadow price-based flow-gate rights (FGR), used to hedge against congestion charges arising on a specific line~\cite{chao:1996}. In either case the rents (obligations) are derived from market equilibrium price differences that arise due to network congestion.

More recently, there has been a growing interest in the development of financial property rights for energy storage systems (ESS)~\cite{taylor:2015,brijs:2016,munoz-alvarez:2017} as a means to optimize their deployment, so as to maximize welfare benefits~\cite{sioshansi:2009}, independently of the ownership structures, and to remunerate the investment. In face of the increasing penetration of RES, the need for increased system flexibility provided by ESS has led many systems to consider ESS as a transmission asset~\cite{pjm:2012,caiso:2018,miso:2018}. Discussion has focused on how to foster efficient investment in ESS and how investments should be remunerated: on cost recovery basis or by merchant revenues in the energy market transactions~\cite{sioshansi:2012}. Motivations for storage rights are analogous to their transmission counterparts, including recovery of investment costs in the market, and the possibility to hedge against LMP volatility over time in the form of forward energy contracts covered by storage capacity.

Several classes of financial property rights for ESS have appeared in the literature. As opposed to their transmission counterparts, for which its impracticality has been shown~\cite{hogan:1992} as a result of the physics of electricity, the feasibility of physical storage rights to ESS was analysed in~\cite{brijs:2016} for market allocation of storage resources in the different energy and ancillary services markets. Two other classes of financial rights to storage analogous to FGRs were introduced in~\cite{taylor:2015}, namely, energy capacity (ECR) and power capacity rights (PCR). In~\cite{munoz-alvarez:2017}, a different class of financial storage rights (FSR) analogous to point-to-point FTRs was alternatively proposed.

The importance of storage in power systems goes back a long time with the employment of hydroelectric power plants, because of both its economic and technical characteristics in support of peak shaving and system flexibility. Moreover, the ability to arbitrage prices in time makes hydroelectricity hold a parallel with modern ESS technologies. However, there are inherent differences between these storage technologies, most notably the fact that hydropower plants have a primary energy source which enables replenishing storage, as opposed to ESS which require electric charging and, in the case of cascaded reservoirs, the hydrological link between storage resources, which gives rise to operational dependency, especially in weekly or monthly timescales.

As power system stability is becoming more complex with the increased variability in generation resulting from higher RES penetration, and, consequently, system flexibility requirements increase, the importance of hydroelectricity flexibility becomes more prevalent, as thoroughly surveyed in the literature~\cite{yang:2018}. Beyond the habitual mechanical and electrical strains inherent to all participating plants, regardless of reservoir size, additional impact is expected on those with smaller storage due to their reduced capacity to arbitrage price differences in a span of several hours or days. Moreover, the profitability of smaller hydroelectric power generation projects has been shown to be mostly dependent on feed-in tariffs~\cite{cunha:2014,redpath:2015}.

\subsection{Contributions and organization}%
\label{sec:introduction:contribution}

This paper proposes a generalized storage model for the theoretical financial rights framework proposed in~\cite{munoz-alvarez:2017}, extending its applicability to hydroelectricity. In that framework, it is shown how ownership and operation of storage can be decoupled with use of FSR and ECR for potential improvements in their valuation and risk management. Our proposed generalization takes into consideration reasonable assumptions regarding the nonlinear productivity and short-term operation properties of hydroelectric power plants, while maintaining convexity of the underlying multi-period economic dispatch model, thus preserving revenue adequacy. The remainder of the text is organized as follows. Section~\ref{sec:hydro} reviews the concepts of hydroelectricity, and presents a rationale for simplification of the production function in the short-term operation. Section~\ref{sec:fsr} summarizes the financial rights framework with a generalized storage model resulting from the assumptions laid out in the preceding section. Finally, Section~\ref{sec:discussion} and Section~\ref{sec:conclusion} discuss practical aspects of financial rights for hydroelectricity, and propose future work, respectively.

%% file: hydro.tex
\section{Hydroelectricity}
\label{sec:hydro}

\begin{figure*}
	\centering
	\subfloat[][Trapezoidal.]{\input{trapezoid.tikz}\label{fig:hydro:reservoir:trapezoid}}
	\quad
	\subfloat[][Cuboidal.]{\input{cuboid.tikz}\label{fig:hydro:reservoir:cuboid}}
	\quad
	\subfloat[][Planar.]{\input{planar.tikz}\label{fig:hydro:reservoir:plane}}
	\caption{Simplified geometric representations of different reservoir models: \protect\subref{fig:hydro:reservoir:trapezoid} in a trapezoidal model the surface area varies with respect to the reservoir height, whereas in the \protect\subref{fig:hydro:reservoir:cuboid} cuboidal model, surface area remains constant with variable height. Conversely, in a \protect\subref{fig:hydro:reservoir:plane} planar model the surface area varies indefinitely with respect to the volume of water, such that height remains constant. In all cases the volume of water is variable and subject to operating limits representing reservoir capacity, either physical or allocated.}%
	\label{fig:hydro:reservoir}
\end{figure*}

Hydroelectricity is obtained from the kinetic energy carried by the movement of water, which is itself harnessed by turbines connected to electrical generators. This kinetic energy can be derived from the gravitational potential energy provided by the construction of river dams that elevate the water upstream of their location, thus forming an upper reservoir that will increase head, and thus the potential energy required for the desired hydroelectric power output. The upper hydraulic head $\varphi(v)$, or forebay height, is a function of the volume of stored water. Immediately downstream of the dam, a lower reservoir is formed by the water released from the upper one, and the resulting tailrace height can be described by a function $\vartheta(q)$. The ideal power output can then be recast as a function of the volume of water in the upper reservoir and discharge:
\begin{displaymath}
	p(v,q) = \rho \cdot g \cdot \left[\varphi(v) - \vartheta(q)\right] \cdot q,
\end{displaymath}
where $\rho$ and $g$ are constants representing the density of water and gravity acceleration, respectively.

It is clear from the equation above that the relationship between $v$ and $q$ will define the mathematical properties of the ideal power output function. Vieira et al.~\cite{vieira:2015} have shown that reasonable assumptions regarding the intrinsic characteristics of river topographies result in strongly increasing and strictly positive $\varphi(v)$ and $\vartheta(q)$ volume-to-height functions. If the area $A$ of a water reservoir surface increases with volume\footnote{\emph{Id est}, $\partial A/\partial v > 0$ in the case of the upper reservoir, or $\partial A/\partial q > 0$ in the case of the lower reservoir.} then its geometry can be approximated by a trapezoid, and the resulting volume-to-height function, by a concave higher order polynomial with positiveness and monotonicity properties, and $\partial^2\varphi/\partial v^2, \, \partial^2 \vartheta/\partial q^2  \leqslant 0$, in the operating range. It is also possible to describe the reservoir models by further simplifications, as illustrated in~\figurename~\ref{fig:hydro:reservoir}. For instance, a cuboidal model results in a linear volume-to-height function where $\partial A/\partial v = 0$, whereas a constant height necessarily implies a planar reservoir model for which $\partial A/\partial v \approx \infty$.

The hydro-electromechanical system defining the hydroelectric power station is subject to energy losses whilst it is converted all the way to electricity. Head, and thus potential energy, is lost due to entrance, pipe bending, and friction losses in the penstock that leads the water from the upper reservoir to the turbine vanes, such that:
\begin{displaymath}
	p(v,q) = \rho \cdot g \cdot h(v,q) \cdot q,
\end{displaymath}
where $h(\cdot) = \left[ \varphi(v) - \vartheta(q) - \varsigma(q) \right]$ denotes the net head, and $\varsigma(q)$ represents a head loss due to penstock losses as a function of water discharge. This water reaches a turbine-generator unit designed for maximum efficiency for given reference values of head and discharge. An efficiency function $\epsilon(v,q)$ can thus be accounted for in the power output equation, such that:
\begin{equation}
	p(v,q) = \pi(v,q) \cdot q%
	\label{eq:hydro}
\end{equation}
where $\pi(.)$ is the hydroelectric productivity, defined as the unit of power output per unit of water discharge:
\begin{displaymath}
	\pi(v,q) = \frac{p(v,q)}{q} = \rho \cdot g \cdot \epsilon(v,q) \cdot h(v,q).
\end{displaymath}

Because the hydroelectric power generation function~\eqref{eq:hydro} is effectively nonlinear and, most notably, non-convex, it is commonly simplified in order to allow for its use in convex optimization algorithms that seek to solve common power systems planning and operation problems in timescales of hours to years. These simplifications will incur errors which may be acceptable if the underlying modeling assumptions are reasonably consistent with the physical characteristics of the system, as well as the timescales under consideration.

In simple terms, as timescales approach real-time, such as in primary frequency control, it becomes reasonable to assume a constant head with respect to both $v$ and $q$, depending on the modeling objectives, whereas efficiency remains variable. On the other hand, as timescales progress from hours or days to several months, with discrete decision time steps measured in terms of a single or even multiple weeks, as in reservoir operation, the hydroelectric power generation function may, conversely, be reasonably represented by variable head and constant efficiency. In both cases productivity is variable.

The most prominent productivity factor in reservoir operation is head. In the long run, optimal reservoir operation for economic dispatch in electric power systems is eventually a problem of deciding on forebay height subject to the seasonal availability of water, and therefore, depending on the reservoir geometry and relative capacity with respect to inflows, a reasonable formulation of~\eqref{eq:hydro} will take variable productivity with respect to forebay height into consideration. In such timescales, factors with short-term variabilities are commonly held constant, e.g. average values. This is the case with turbine-generator efficiency~\cite{martins:2014}, or even tailrace height~\cite{vieira:2015} if the decision timeframes are long enough, e.g. weekly or monthly time steps. If uncertainties regarding water inflows are explicitly considered in the model, simplifications of~\eqref{eq:hydro} reported in the literature go as far as fixing productivity to accommodate for stochastic optimisation~\cite{pereira:1991}.

On the other hand, in the short run, productivity is more sensibly influenced by short-term variables. Depending on reservoir geometry, its relative capacity with respect to water discharge, and turbine type, head and efficiency factors might be more prominently influenced by tailrace height, as forebay height changes slowly in timescales of hours or days, especially and more significantly in very large upper reservoirs. In this case, holding forebay height constant in~\eqref{eq:hydro} is a common approach~\cite{lund:1999,paredes:2015,dal-santo:2016}, although the modeling of head variability by forebay height is not uncommon~\cite{diniz:2008,catalao:2009}.

\subsection{Power generation model}%
\label{sec:hydro:generation}
The nonlinearity of the hydroelectric power generation model arises from the variable productivity commonly formulated as a polynomial function $\pi:\mathcal{V}\times\mathcal{Q}\rightarrow\Real_{+}$, i.e. a function of head and water discharge. Its order can be reduced in two main ways, depending on the timescales involved, such that $\tilde{\pi}:\mathcal{V}\rightarrow\Real_+$ or $\tilde{\pi}:\mathcal{Q}\rightarrow\Real_+$, whether long- or short-term operation is concerned, respectively. In the case of day-ahead markets a model of $\tilde{\pi}(q)$ that retains convexity is desirable. This can be achieved by reducing~\eqref{eq:hydro} down to a second degree polynomial with the following simplifying assumptions: a planar upper reservoir model, i.e. constant forebay height, cuboid geometry of the tailrace canal, and constant efficiency. Thus:
\begin{IEEEeqnarray*}{rcl}
	p(q) & \; = \; & \tilde{\pi}(q) \cdot q,\\
	\tilde{\pi}(q) & \; = \; & \kappa \cdot \tilde{h}(q),\\
	\tilde{h}(q) & \; = \; & \tilde{\varphi} - \varphi(q) - \tilde{\varsigma},\\
	\kappa & \; = \; & \rho \cdot g \cdot \epsilon(\tilde{v}, \tilde{q}),
\end{IEEEeqnarray*}
where $\kappa$ denotes constant efficiency factor for fixed values of $v$ and $q$, and $\tilde{\varphi} > \tilde{\varsigma} \geqslant 0$ denote given forebay height and penstock head loss values, respectively. Moreover, tailrace height varies linearly with respect to $q$, thus:
\begin{equation}
	\label{eq:hydro:generation:theta}%
	\vartheta(q) = \theta_0 + \theta_1 \cdot q,
\end{equation}
where $\tilde{\varphi} > \theta_0 > \tilde{\varsigma} \geqslant 0$, such that $\tilde{h}:\mathcal{Q}\rightarrow\Real_+$. In the case of of impulse turbines, $\theta_1=0$ and, consequently, tailrace height remains constant, and thus $p(q)$ is linear. On the other hand, in the case of reaction turbines, $\theta_1>0$ and $p(q)$ is quadratic, assuming the following form:
\begin{equation}
	\label{eq:hydro:generation:p}
	p(q) = \kappa \cdot \left[ \left(\tilde{\varphi} - \theta_0 - \tilde{\varsigma}\right) - \theta_1 \cdot q \right] \cdot q,
\end{equation}
where $p(q)$ is positive, continuous, twice differentiable, and increasing for $q\in\mathcal{Q}\subset\Real_+$, such that $p(0)=0$. Since we are interested in the consideration of hydroelectric power as a storage resource, we shall alternatively denote $p(q)$ simply as $u(q)$, or hydroelectric storage power discharge.

Let $\alpha = \kappa \cdot \theta_1$, and $\beta = \kappa \cdot \left(\tilde{\varphi} - \theta_0 - \tilde{\varsigma}\right)$, such that $\alpha,\beta\geqslant 0$. We shall rewrite the hydroelectric power generation function~\eqref{eq:hydro:generation:p} as:
\begin{equation}
	\label{eq:hydro:generation:u}%
	u(q) = -\alpha \cdot q^2 + \beta \cdot q.
\end{equation}
Inversely, we can write $q_i(u_i)$ for a given plant $i$, $0 \leqslant \vetor{u}[i] \leqslant \maximum{\vetor{u}}[i]$, in a change of variables, with water discharge as a function of power discharge in a convex quadratic Taylor series approximation at $u_i=0$, such that $q_i(0) = 0$:
\begin{equation}
	\label{eq:hydro:generation:q}%
	q_i(u_i) = a_i \cdot u_i^2 + b_i\cdot u_i + \cancelto{0}{\frac{\beta_i - \sqrt{\beta_i^2}}{2\cdot\alpha_i}},
\end{equation}
where $a_i = \alpha_i\cdot\sqrt{\beta_i^2}/\beta_i^4$, and $b_i = 1/\sqrt{\beta_i^2}$.

\subsection{Storage model}%
\label{sec:hydro:storage}
Most of the dissimilarities between ESS and hydroelectricity regarding the mathematical formulation concern the storage model. In the former, the energy storage is directly charged from the grid. In the latter, because energy is stored in its primary source--water, the reservoir is exogenously (with respect to the grid) replenished from the river basin streamflows. These inflows can be either natural or controllable from upstream reservoirs, thus establishing either weakly or strongly connected reservoirs. In weakly connected reservoirs the routing effect is sufficiently slow not to be captured in the time window of the economic dispatch problem, whereas between strongly connected reservoirs this relationship must be considered for proper storage modeling. In this latter case, in the timescales of interest, the underlying dynamics is commonly linearly approximated in the discrete time periods.

In its most simple formulation, the energy storage model of ESS is represented by steady-state linear, lossless dynamics:
\begin{equation}
	\label{eq:fsr:storage:dynamic:z}%
	Z_{i,t+1} = Z_{i,t} -  U_{i,t},
\end{equation}
where $Z_{i,0}$ is given and represents the initial storage condition, whose reformulation in terms of $U_{i,\cdot}$ subject to minimum and maximum storage capacities $\minimum{Z},\maximum{Z}\in\Real^{n\times T}_{+}$, respectively:
\begin{displaymath}
	\minimum{Z}_{i,t} \leqslant Z_{i,t} = Z_{i,0} - \sum_{\tau=1}^t U_{i,\tau} \leqslant \maximum{Z}_{i,t},
\end{displaymath}
gives rise to a lower triangular matrix $L\in\Real^{T\times T}$, whose nonzeros equal $-1$, such that $\minimum{\vetor{z}}[i] \leqslant L \vetor{u}[i] \leqslant \maximum{\vetor{z}}[i]$.

On the other hand, representing either isolated or weakly connected reservoirs, the hydroelectric storage model can be formulated as follows, while accounting for the power output to water discharge conversion function:
\begin{equation}
    \label{eq:discussion:storage}%
    \minimum{\vetor{z}}[i] \leqslant L\vetor{q}[i] + \sum_{\forall j\in\Omega}\vetor{q}[j] + \vetor{y}[i] \leqslant \maximum{\vetor{z}}[i],
\end{equation}
where $\vetor{y}[i]\in\Real^T$ denotes given incremental inflows to reservoir $i$ over time, $\minimum{\vetor{z}}[i]$ and $\maximum{\vetor{z}}[i]$ are given as minimum and maximum reservoir volume values, respectively, $\Omega$ represents the set of all reservoirs $j$ immediately upstream of $i$, and $\vetor{q}[i]:\Upsilon_i\subset\Real^{T}\rightarrow\Real^{T}$ represents the quadratic conversion rate function~\eqref{eq:hydro:generation:q}, as follows:
\begin{equation}
	\label{eq:hydro:storage:q}%
    \vetor{q}[i] = \transpose{\left\{\transpose{\vetor{u}[i]} \left( A_i \text{diag}(\vetor{u}[i]) + B_i \right)\right\}},
\end{equation}
where $A_i,B_i\in\Real^{T\times T}_+$ are given by $a_i I$ and $b_i I$, respectively, and $I$ denotes the identity matrix of suitable dimensions.

%% file: trapezoid.tikz
\begin{tikzpicture}[scale=.75, z={(-.75, -.25)}]
	\node [draw=none] (0) at (-2.5, -2) {};
	\node [draw=none] (1) at (-1.5, -4) {};
	\node [draw=none] (2) at (0.5, -4) {};
	\node [draw=none] (3) at (1.5, -2) {};
	\node [draw=none] (4) at (0, 0) {};
	\node [draw=none] (5) at (4, 0) {};
	\node [draw=none] (6) at (3, -2) {};
	\node [draw=none] (7) at (1, -2) {};
	\draw (0.center) to (1.center);
	\draw (1.center) to (2.center);
	\draw (2.center) to (3.center);
	\draw (0.center) to (3.center);
	\draw (0.center) to (4.center);
	\draw (3.center) to (5.center);
	\draw (4.center) to (5.center);
	\draw (2.center) to (6.center);
	\draw (6.center) to (5.center);
	\draw [style=dashed, color=gray] (4.center) to (7.center);
	\draw [style=dashed, color=gray] (1.center) to (7.center);
	\draw [style=dashed, color=gray] (6.center) to (3.center);
\end{tikzpicture}

%% file: cuboid.tikz
\begin{tikzpicture}[scale=.75, z={(-.75, -.25)}]
	\node [draw=none] (0) at (-2.5, -2) {};
	\node [draw=none] (1) at (-2.5, -4) {};
	\node [draw=none] (2) at (1.5, -4) {};
	\node [draw=none] (3) at (1.5, -2) {};
	\node [draw=none] (4) at (0, 0) {};
	\node [draw=none] (5) at (4, 0) {};
	\node [draw=none] (6) at (4, -2) {};
	\node [draw=none] (7) at (0, -2) {};
	\node [draw=none] (8) at (2.75, -2) {};
	\draw (0.center) to (1.center);
	\draw (1.center) to (2.center);
	\draw (2.center) to (3.center);
	\draw (0.center) to (3.center);
	\draw (0.center) to (4.center);
	\draw (3.center) to (5.center);
	\draw (4.center) to (5.center);
	\draw (2.center) to (6.center);
	\draw (6.center) to (5.center);
	\draw [style=dashed, color=gray] (4.center) to (7.center);
	\draw [style=dashed, color=gray] (1.center) to (7.center);
	\draw [style=dashed, color=gray] (6.center) to (3.center);
\end{tikzpicture}

%% file: planar.tikz
\begin{tikzpicture}[scale=.75, z={(-.75, -.25)}]
	\node [draw=none] (0) at (-2.5, -2) {};
	\node [draw=none] (3) at (1.5, -2) {};
	\node [draw=none] (4) at (0, 0) {};
	\node [draw=none] (5) at (4, 0) {};
	\draw (0.center) to (3.center);
	\draw [style=dashed, color=gray] (0.center) to (4.center);
	\draw [style=dashed, color=gray] (3.center) to (5.center);
	\draw (4.center) to (5.center);
\end{tikzpicture}

%% file: fsr.tex
\section{Financial Rights}%
\label{sec:fsr} 

The financial rights framework introduced in~\cite{munoz-alvarez:2017} is summarized in this section, along with our proposed generalized convex quadratic storage model, and a brief discussion about its implications to the theoretical results. Some changes in notation were required in order to accommodate for our hydroelectricity model. In this framework, both transmission and storage are centrally operated and available as open access resources such that welfare maximization is sought after subject to constraints on the use of transmission and storage.

Let the transmission network be represented by a connected, directed graph with $n$ buses and $m$ lines over a number of $T$ discrete time periods, then the multi-period energy-constrained economic dispatch follows:

\begin{IEEEeqnarray}{rcl}
	\IEEEyesnumber \label{eq:fsr:mped} \IEEEyessubnumber*
	\min_{P, U} \quad \quad \sum_{t=1}^{T} \sum_{i=1}^n c_i\left(P_{i,t}, t\right) \label{eq:fsr:mped:cost}\\
	\text{subject to} \qquad \mathbf{1}^\intercal \left\{\vetor{p}[t]+\vetor{u}[t]\right\} & \; = \; & 0, \quad t=1,\ldots, T \label{eq:fsr:mped:power}\\
	  G\left\{\vetor{p}[t] + \vetor{u}[t]\right\} & \; \leqslant \; & \vetor{d}, \quad t=1,\ldots, T\label{eq:fsr:mped:flow}\\
	  \minimum{\vetor{z}}[i] \leqslant L\vetor{q}[i] + \sum_{\forall j \in \Omega} \vetor{q}[j] + \vetor{y}[i]& \; \leqslant \; & \maximum{\vetor{z}}[i], \, i=1,\ldots, n\label{eq:fsr:mped:storage}
\end{IEEEeqnarray}
where $P\in\Psi\subset\Real^{n\times T}$ and $U\in\Upsilon\subset\Real^{n\times T}$ define the multi-period economic dispatch problem variables, respectively representing power generation (consumption), if $P_{i,t} > 0$ $\left(P_{i,t} < 0\right)$, and storage power discharge (charge), if $U_{i,t} > 0$ ($U_{i,t} < 0$), at a bus $i$ and period $t$, such that $\vetor{p}[i],\vetor{u}[i]\in\Real^{T}$ define the vectors associated with the respective $i$-th rows, and $\vetor{p}[t],\vetor{u}[t]\in\Real^{n}$, with $t$-th columns.

\begin{remark}
	Sets $\Psi$ and $\Upsilon$ define implicit power generation (consumption) and storage (dis)charge operating ranges.
\end{remark}

In the objective function~\eqref{eq:fsr:mped:cost} of the above problem, offers and bids are represented by increasing, differentiable, and convex functions $c_i(\cdot)>0$ for $P_{i,t}>0$, and concave functions $-c_i(\cdot)<0$ for $P_{i,t}<0$, respectively. Net power injections at all buses are subject to constraints on a lossless linear DC power approximation of the transmission network, i.e. power balance equations~\eqref{eq:fsr:mped:power}, and maximum power flows~\eqref{eq:fsr:mped:flow} in every line, where $G\in\Real^{2m\times n}$ is the matrix of shift factors, and $\vetor{d}\in\Real^{2m}_+$ is a vector representing directed maximum line capacities. If $A_i = B_i = 0\in\Real^{T\times T}$, then~\eqref{eq:fsr:mped:storage} reduces to the specific ESS case where $\minimum{\vetor{z}}[i] \leqslant L \vetor{u}[i] \leqslant \maximum{\vetor{z}}[i]$.

In efficient market equilibrium there exists a set of locational marginal prices $\Lambda\in\Real^{n\times T}$ derived from the first-order optimality conditions of the Lagrangian function:
\begin{align*}
	\lagrangian(\cdot) =&  \sum_{t=1}^{T} \sum_{i=1}^n c_i\left(P_{i,t}, t\right) - \gamma_t \transpose{\ones} \left\{\vetor{p}[t] + \vetor{u}[t]\right\} \\
		& \quad - \transpose{\vetor{\mu}[t]} G \left\{\vetor{p}[t] + \vetor{u}[t] - \vetor{d} \right\} \\
		& \quad + \sum_{i=1}^n \transpose{\minimum{\vetor{\eta}}[i]} \{\minimum{\vetor{z}}[i] - L\vetor{q}[i] - \sum_{\forall j \in \Omega} \vetor{q}[j] - \vetor{y}[i]\} \\
		& \qquad + \transpose{\maximum{\vetor{\eta}}[i]} \{ L\vetor{q}[i] + \sum_{\forall j \in \Omega} \vetor{q}[j] + \vetor{y}[i] - \maximum{\vetor{z}}[i] \}, \numberthis \label{eq:fsr:lagrangian}
\end{align*}
such that $\vetor{\lambda}[t] = \gamma_t \ones - G\vetor{\mu}[t]$.

If the optimal solution occurs with either transmission or storage congestion in~\eqref{eq:fsr:mped:flow} or~\eqref{eq:fsr:mped:storage}, respectively, the system operator will collect a positive merchandising surplus (zero otherwise), defined as $-$trace$(\transpose{\Lambda}P)$, or equivalently:
\begin{align*}
	\sum_{t=1}^{T} \transpose{\vetor{\mu}[t]} G \left\{ \vetor{p}[t] + \vetor{u}[t]\right\}\\
	+ \sum_{i=1}^{n} \transpose{\left( \maximum{\vetor{\eta}}[i] - \minimum{\vetor{\eta}}[i] \right)} (L \vetor{q}[i] + \sum_{\forall j\in\Omega} \vetor{q}[j] + \vetor{y}[i]) &\geqslant 0, \numberthis \label{eq:fsr:ms}
\end{align*}
since the convexity of~\eqref{eq:hydro:storage:q} holds the above inequality true.

The revenue generated by a merchandising surplus is distributed to holders of financial property rights, as is the practice in several US electricity markets under congested transmission, in the form of either a FTR or a FGR. These financial rights are assigned to their holders by means of regular auctions.
\begin{definition}
	Let $l$ denote a directed line, and $i$, $k$ represent an injection and withdrawal bus, respectively. A point-to-point FTR is a right (obligation) for its holder to a rent (liability) equal to $\transpose{\left( \vetor{\lambda}[k] - \vetor{\lambda}[i] \right)} \vetor{r}_{ik}$, where $\vetor{r}_{ik}\in\Real^{T}$ is a power profile of a symmetrical injection and withdrawal, whereas a FGR is a right for its holder to a rent equal to $\transpose{\vetor{\mu}[l]}\vetor{f}[l]$, where $\vetor{f}[l]\in\Real^{T}$ is a power profile
\end{definition}

In the case of financial property rights to storage, they can be defined as either FSR~\cite{munoz-alvarez:2017} or ECR~\cite{taylor:2015}.
\begin{definition}
	Let $i$ and $k$ denote a target storage asset and a withdrawal bus, respectively. A time-to-time FSR is a right (obligation) for its holder to a rent (liability) equal to $\transpose{\vetor{\lambda}[k]} \vetor{s}[k]$, where $\vetor{s}[k]\in\Real^{T}$ is a power profile, whereas an ECR is a right for its holder to a rent equal to $\transpose{\maximum{\vetor{\eta}}[i]} \vetor{e}[i]$, where $\vetor{e}[i]\in\Real^{T}_+$ denotes a power profile.
\end{definition}
\begin{definition}
    Let $(\vetor{R}_{ik},\vetor{F}[l],\vetor{S}[i],\vetor{E}[i])$ denote a collection of financial transmission and storage rights defined by the sum of all issued rights---e.g. by means of auctions, $\vetor{r}_{ik}$, $\vetor{f}[l]$, $\vetor{s}[i]$, and $\vetor{e}[i]$, respectively.
\end{definition}
It can be demonstrated under the generalized storage model that, analogously to~\cite{munoz-alvarez:2017}, and upon optimality, the rents due to a collection of financial rights $(\vetor{R}_{ik},\vetor{F}[l],\vetor{S}[i],\vetor{E}[i])$ are covered by the merchandising surplus, provided the convexity of~\eqref{eq:fsr:mped}---thus guaranteeing system operator revenue adequacy:
\begin{align*}
	-\sum_{t=1}^T \transpose{\vetor{\lambda}[t]}\vetor{p}[t] \geqslant & \sum_{i,k}^n \transpose{\left( \vetor{\lambda}[k] - \vetor{\lambda}[i] \right)} \vetor{R}_{ik} + \sum_{l=1}^{2m} \transpose{\vetor{\mu}[l]} \vetor{F}[l] \\
	 & \quad + \sum_{i=1}^n \transpose{\vetor{\lambda}[i]} \vetor{S}[i] + \transpose{\maximum{\vetor{\eta}}[i]} \vetor{E}[i] \numberthis \label{eq:fsr:ra}
\end{align*}
as long as $(\vetor{R}_{ik},\vetor{F}[l],\vetor{S}[i],\vetor{E}[i])$ solves the generalized simultaneous feasibility test:
\begin{IEEEeqnarray}{rcl}
    \IEEEyesnumber \label{eq:fsr:sft} \IEEEyessubnumber*
    \mathbf{1}^\intercal \left\{\vetor{R}[t] - \vetor{S}[t] + \vetor{u}[t]\right\} & \; = \; & 0, \qquad\qquad\; \forall t \label{eq:fsr:sft:power}\\
	  G\left\{\vetor{R}[t] - \vetor{S}[t] + \vetor{u}[t]\right\} & \; \leqslant \; & \vetor{d} - \vetor{F}[t], \quad\, \forall t\label{eq:fsr:sft:flow}\\
	  \minimum{\vetor{z}}[i] \leqslant L\vetor{q}[i] + \sum_{\forall j\in\Omega}\vetor{q}[j] + \vetor{y}[i]& \; \leqslant \; & \maximum{\vetor{z}}[i] - \vetor{E}[i], \; \forall i\label{eq:fsr:sft:storage}
\end{IEEEeqnarray}

%% file: discussion.tex
\section{Discussion}%
\label{sec:discussion}

In the financial rights framework presented, a participating hydroelectric power plant would be centrally dispatched as an open access resource by the system operator to its full available capacity, making itself available to the energy and ancillary services markets so as to maximize its economic value to the system. Its practical implementation in the day-ahead markets could, however, present a few adaptations to account for characteristics that set hydroelectric power apart from other energy storage technologies.

In this setting financial storage rights would be offered in combination with energy bids, in which the hydroelectric generator would bid the minimum price at which it is willing to supply a specified amount of energy on the following day, whether to its full capacity or not. If the bid is accepted, the FSR would then allow the system operator to determine how the hydroelectric power plant is dispatched in all periods of the day-ahead market to maximize system value, either providing energy or ancillary services.

This arrangement would substantially simplify the bid-making process of the hydropower generators, since they would no longer have to submit bids for each time interval of the day-ahead market based on their projections of the hourly load and of the supply from the remaining generators. The FSR bids would only require a single price-quantity bid for the entire day. Thus, these generators would only need to determine the opportunity costs of generating a specified amount of energy in the day-ahead market compared against their long-term operation planning.

Alternatively, these bids could also include capacity constraints for these hydroelectric power plants, as well as production function parameters, so as to enable the system operator to take into consideration its productivity differences based on how the reservoirs are operated. An initial dispatch order would be determined considering flat bids for these generators, i.e. a constant supply value during all time periods of the day-ahead market. The marginal prices for each of these periods would then be determined given the initial dispatch order.

Next, hydroelectric power generation would be reallocated from periods with lower marginal prices to those with higher marginal prices. The value of FSR could then correspond to the cost savings provided---the difference between the marginal prices of the initial dispatch order, and the resulting marginal costs after the energy reallocation.

%% file: conclusion.tex
\section{Conclusion}%
\label{sec:conclusion}

A generalized storage model for the financial rights framework proposed in~\cite{munoz-alvarez:2017}, extending its applicability to hydroelectricity, was presented. The model accounts for variable productivity resulting from reasonable assumptions regarding nonlinear short-term production functions, while enabling the representation of hydroelectric power as a series of linear injections, resulting in a convex quadratic multi-period economic dispatch problem that preserves its revenue adequacy properties. A brief discussion about the practical aspects for the implementation of FSR for hydroelectricity was also presented. Future work includes numerical case study simulations for the analysis of the power generation model approximations, as well as the applicability patterns of these financial rights for market power mitigation in cascaded reservoirs.